# Summation rules for the values of the Riemann zeta–function and generalized harmonic numbers obtained using Laurent developments of polygamma functions and their products

Sergey K. Sekatskii

*Laboratory of Biological Electron Microscopy, IPHYS, BSP 419, Ecole Polytechnique Fédérale de Lausanne, and Dept. of Fundamental Biology, Faculty of Biology and Medicine, University of Lausanne, CH1015 Lausanne-Dorigny, Switzerland.*
*E-mail :* serguei.sekatski@epfl.ch

Following the Mellin and inverse Mellin transforms technique presented in our paper arXiv:1606.02150 (NT), we have established close forms of Laurent series expansions of products of bi- and trigamma functions $\psi'(z)\psi(-z)$ and $\psi'(z)\psi'(-z)$. These series were used to find summation rules which include generalized harmonic numbers of first, second and third powers and values of the Riemann zeta-functions at integers / Bernoulli numbers, for example $2\sum_{k=1}^{\infty}\frac{H_k^{(2)}}{k^3}=6\varsigma(2)\varsigma(3)-9\varsigma(5)$. Some of these rules were tested numerically.

## Introduction

In a recent paper [1] we have established numerous integral representations of the Riemann zeta-function and other related functions as Mellin transforms of appropriate trigonometrical functions. Then, using inverse Mellin transform technique, we obtained Laurent expansions of the trigonometrical functions involved, which in their turn were used to establish summation rules for the Riemann zeta-functions and closely related with them Bernoulli, Euler, or some related numbers. In the last Section of [1], not trigonometrical functions but digamma functions and its derivatives were used for the same purposes, and we sketched the program of further applications of this approach. In this short Note we report the first realizations of such a program.

## 1. Summation rules

*1.1 Product of trigamma functions*

Let us start with the function $\psi'(z)\psi'(-z)$ which is an "even product" of two trigamma functions (derivatives of the digamma functions, see e.g. [2] for necessary definitions and basic properties). We know the Laurent series

$$\psi(x) = -\frac{1}{x} - \gamma + \sum_{n=1}^{\infty}(-1)^{n-1}\varsigma(n+1)x^n \qquad (1),$$

valid for $|x|<1$ (see again [2]; these series were also obtained in [1]). Here $\varsigma(s)$ is the Riemann zeta-function (see e.g. [3] for its main properties) and $\gamma$ is Euler – Mascheroni constant. From (1) we have

$$\psi'(x) = \frac{1}{x^2} + \sum_{n=0}^{\infty}(-1)^n(n+1)\varsigma(n+2)x^n \qquad (2),$$

and it is easy to see that the function $F(z) := \psi'(z)\psi'(-z) - \frac{1}{z^4} - \frac{2\varsigma(2)}{z^2} - \varsigma^2(2) - 6\varsigma(4)$ is *$O(z^2)$* at zero and *O(1)* at infinity (trigamma function is *O(1/z)* at infinity [2]), hence the contour integral $\int_C z^{-s-1}(\psi'(z)\psi'(-z) - \frac{1}{z^4} - \frac{2\varsigma(2)}{z^2} - \varsigma^2(2) - 6\varsigma(4))dz$ exists in the strip $0<\mathrm{Re}s<2$. Here the contour $C$ consists of the line $-iX$, $+iX$ and demi-circle connecting these points and lying in the right complex semi-plane; $X$ is real and tends to infinity. (If a pole of the function under integral sign occurs on the contour, just change the value of $X$ a bit to avois this). In this strip, this integral is equal to (for consistence, we repeat the material from [1]):

$$\int_C z^{-s-1}\left(\psi'(z)\psi'(-z) - \frac{1}{z^4} - \frac{2\varsigma(2)}{z^2} - \varsigma^2(2) - 6\varsigma(4)\right)dz$$

$$= -i(i^{-s-1} + (-i)^{-s-1})\int_0^{\infty} x^{-s-1}\left(\psi'(ix)\psi'(-ix) - \frac{1}{x^4} + \frac{2\varsigma(2)}{x^2} - \varsigma^2(2) - 6\varsigma(4)\right)dx$$

$$= -2i\cos(\frac{\pi(s+1)}{2})\int_0^\infty x^{-s-1}\left(\psi'(ix)\psi'(-ix) - \frac{1}{x^4} + \frac{2\varsigma(2)}{x^2} - \varsigma^2(2) - 6\varsigma(4)\right)dx$$

$$= 2i\sin(\frac{\pi s}{2})\int_0^\infty x^{-s-1}\left(\psi'(ix)\psi'(-ix) - \frac{1}{x^4} + \frac{2\varsigma(2)}{x^2} - \varsigma^2(2) - 6\varsigma(4)\right)dx \qquad (3).$$

From the other side, this contour integral can be calculated using Cauchy residues theorem as equal to $2\pi i(\sum_{n=1}^{\infty} n^{-s-2}\psi'(n)(s+1) - \sum_{n=1}^{\infty}\psi''(n)n^{-s-1})$: in the vicinity of points *z*=1, 2, 3… the function $\psi(-z)$ is given by [1, 2]

$$\psi(-n-x) = \frac{1}{x} + H_n - \gamma + \sum_{k=1}^{\infty}((-1)^k H_n^{(k+1)} - \varsigma(k+1))x^k \qquad (4),$$

whence its derivative has here poles of the second order having Laurent expansion $-\frac{1}{x^2} + O(1)$. Here $H_n^{(k)} = \sum_{l=1}^{n}\frac{1}{l^k}$ is generalized harmonic number, $H_0^{(k)} = 0$ by definition. Below we will omit superscript 1 for $H_n^{(1)} = \sum_{l=1}^{n}\frac{1}{l}$ writing it simply as $H_n$.

Now we need to introduce the functions

$$\chi_k(s) := \sum_{n=1}^{\infty}\psi^{(k)}(n)n^{-s} \qquad (5),$$

where $\psi^{(k)} = \frac{d^k\psi}{dz^k}$ (similar functions were introduced and considered earlier, see e.g. [4]) to write in a compact form

$$\int_0^\infty x^{-s-1}\left(\psi'(ix)\psi'(-ix) - \frac{1}{x^4} + \frac{2\varsigma(2)}{x^2} - \varsigma^2(2) - 6\varsigma(4)\right)dx = \frac{\pi}{\sin(\frac{\pi s}{2})}((s+1)\chi_1(s+2) - \chi_2(s+1)) \qquad (6).$$

This is a Mellin transform. Inverse Mellin transform $F(ix) = \frac{1}{2\pi i}\int_{c-i\infty}^{c+i\infty}\pi x^s\left(\frac{(s+1)\chi_1(s+2) - \chi_2(s+1)}{\sin(s\pi/2)}\right)ds$ (it is easy to see, that all conditions to apply the corresponding theorem hold, see [1], [5], we take *0<c<2*), by virtue of Cauchy residues theorem (simple poles are located at *s=2n*, *n*=1, 2, 3… with the residues of $2/\pi$ and alternating signs) then easily supplies the following Laurent series development valid for $|x|<1$:

$$\psi'(iz)\psi'(-iz) - \frac{1}{z^4} + \frac{2\varsigma(2)}{z^2} - \varsigma^2(2) - 6\varsigma(4) = 2\sum_{n=1}^{\infty}(-1)^n((2n+1)\chi_1(2n+2) - \chi_2(2n+1))z^{2n}$$ and thus

$$\psi'(z)\psi'(-z) - \frac{1}{z^4} - \frac{2\varsigma(2)}{z^2} - \varsigma^2(2) - 6\varsigma(4) = 2\sum_{n=1}^{\infty}((2n+1)\chi_1(2n+2) - \chi_2(2n+1))z^{2n} \qquad (7).$$

To obtain the summation rules, this should be compared with the product of $\psi'(x)$, eq. (2), and $\psi'(-x)=\frac{1}{x^2}+\sum_{n=0}^{\infty}(n+1)\varsigma(n+2)x^n$. This gives for the coefficients standing in front of $x^{2n}$, $n$=1, 2, 3…:

$$2(2n+3)\varsigma(2n+4)+\sum_{l=0}^{2n}(-1)^l(l+1)(2n-l+1)\varsigma(l+2)\varsigma(2n-l+2)=2(2n+1)\chi_1(2n+2)-2\chi_2(2n+1)$$

(8).

Now, we remind that $\psi'(z)=\sum_{l=0}^{\infty}\left(\frac{1}{l+z}\right)^2$ and $\psi''(z)=-2\sum_{l=0}^{\infty}\left(\frac{1}{l+z}\right)^3$. Using harmonic numbers notation, we can write $\psi'(n)=\varsigma(2)-H_{n-1}^{(2)}=\varsigma(2)-H_n^{(2)}+\frac{1}{n^2}$, $\psi''(n)=-2\varsigma(3)+2H_{n-1}^{(3)}=-2\varsigma(3)+2H_n^{(3)}-\frac{2}{n^3}$ and obtain $\chi_1(s)=\varsigma(2)\varsigma(s)-\sum_{n=1}^{\infty}\frac{H_n^{(2)}}{n^s}+\varsigma(s+2)$ and $\chi_2(s)=-2(\varsigma(3)\varsigma(s)-\sum_{n=1}^{\infty}\frac{H_n^{(3)}}{n^s}+\varsigma(s+3))$, that is

$$\chi_1(2n+2)=\varsigma(2)\varsigma(2n+2)-\sum_{n=1}^{\infty}\frac{H_k^{(2)}}{k^{2n+2}}+\varsigma(2n+4)$$

$\chi_2(2n+1)=-2\varsigma(3)\varsigma(2n+1)+2\sum_{k=1}^{\infty}\frac{H_k^{(3)}}{k^{2n+1}}-2\varsigma(2n+4))$. Thus we can recast eq. (8) as the following Proposition:

**Proposition 1.** For *n=1, 2, 3…*, we have

$$\sum_{l=1}^{2n-1}(-1)^l(l+1)(2n-l+1)\varsigma(l+2)\varsigma(2n-l+2)=4\varsigma(3)\varsigma(2n+1)-2(2n+1)\sum_{k=1}^{\infty}\frac{H_k^{(2)}}{k^{2n+2}}-4\sum_{k=1}^{\infty}\frac{H_k^{(3)}}{k^{2n+1}} \quad (10).$$

For *n=2, 3…* this can be written as:

$$\sum_{l=2}^{2n-2}(-1)^l(l+1)(2n-l+1)\varsigma(l+2)\varsigma(2n-l+2)=4(2n+1)\varsigma(3)\varsigma(2n+1)-2(2n+1)\sum_{k=1}^{\infty}\frac{H_k^{(2)}}{k^{2n+2}}-4\sum_{k=1}^{\infty}\frac{H_k^{(3)}}{k^{2n+1}}$$

(10a).

The first two examples of (10), corresponding to respectively $n$=1 and $n$=2, are

$$3\sum_{k=1}^{\infty}\frac{H_k^{(2)}}{k^4}+2\sum_{k=1}^{\infty}\frac{H_k^{(3)}}{k^3}=4\varsigma^2(3) \quad (11)$$

and

$$10\sum_{k=1}^{\infty}\frac{H_k^{(2)}}{k^6}+4\sum_{k=1}^{\infty}\frac{H_k^{(3)}}{k^5}=20\varsigma(3)\varsigma(5)-9\varsigma^2(4) \quad (12).$$

Equalities (11), (12) have been tested numerically.

*1.2 Function* $\psi(x)\psi'(-x)$

Next let us analyse the case of $\psi(x)\psi'(-x)\pm\psi'(x)\psi(-x)$. We trivially have, from earlier established [1] $\psi(x)\psi(-x)=-\frac{1}{x^2}+\gamma^2+\frac{\pi^3}{3}-2\sum_{n=1}^{\infty}\chi_0(2n+1)x^{2n}$ that

$$\psi'(x)\psi(-x)-\psi(x)\psi'(-x)=\frac{2}{x^3}-4\sum_{n=1}^{\infty}n\chi_0(2n+1)x^{2n-1} \qquad (13),$$

here $\chi_0(s):=\sum_{n=1}^{\infty}\psi(n)n^{-s}$ .

To find Laurent series for $\psi(x)\psi'(-x)+\psi'(x)\psi(-x)$, let us introduce the function $F_1(x):=\psi(x)\psi'(-x)+\psi(-x)\psi'(x)+\frac{2\gamma}{x^2}+2\gamma\varsigma(2)+6\varsigma(3)$. We have that in the vicinity of zero $F_1(x)=O(x^2)$, so that the contour integral $\int_C z^{-s-1}F_1(z)dz$ exists in the strip *0*<Re*s*<*2*. Using the same technique as above, we get

$\int_0^{\infty}x^{-s-1}F_1(ix)dx=\frac{\pi}{\sin(\pi s/2)}\sum_{n=1}^{\infty}(2\psi'(n)n^{-s-1}-(s+1)\psi(n)n^{-s+2})$, and then using again the inverse Mellin transform theorem, we obtain Taylor expansion

$$\psi'(ix)\psi(-ix)+\psi(ix)\psi'(-ix)=\frac{2\gamma}{x^2}-2\gamma\varsigma(2)-6\varsigma(3)+2\sum_{n=1}^{\infty}(-1)^n(2\chi_1(2n+1)-(2n+1)\chi_0(2n+2))x^{2n}$$

and thus

$$\psi'(x)\psi(-x)+\psi(x)\psi'(-x)=-\frac{2\gamma}{x^2}-2\gamma\varsigma(2)-6\varsigma(3)-2\sum_{n=1}^{\infty}(2\chi_1(2n+1)-(2n+1)\chi_0(2n+2))x^{2n}$$

(14).

Together with (13), we have

$$\psi'(x)\psi(-x)=\frac{1}{x^3}-\frac{\gamma}{x^2}-\gamma\varsigma(2)-3\varsigma(3)-\sum_{n=1}^{\infty}(2\chi_1(2n+1)-(2n+1)\chi_0(2n+2))x^{2n}-\sum_{n=1}^{\infty}2n\chi_0(2n+1)x^{2n-1}$$

(15)

and

$$\psi(x)\psi'(-x)=-\frac{1}{x^3}-\frac{\gamma}{x^2}-\gamma\varsigma(2)-3\varsigma(3)-\sum_{n=1}^{\infty}(2\chi_1(2n+1)-(2n+1)\chi_0(2n+2))x^{2n}$$

$$+\sum_{n=1}^{\infty}2n\chi_0(2n+1)x^{2n-1} \qquad (16).$$

Summation rules are obtained comparing (15) with the product of (1) and (2), the latter is taken as $\psi'(-x)=\frac{1}{x^2}+\sum_{n=0}^{\infty}(n+1)\varsigma(n+2)x^n$ .

To simplify the logistics, we put here $\varsigma(1)=\gamma$ (formally, of course, but it is known that this is Cauchy principal value $\gamma=\frac{1}{2}\lim_{\varepsilon\to 0}(\varsigma(1+\varepsilon)+\varsigma(1-\varepsilon))$ [3]), and rewrite (1) as $\psi(x)=-\frac{1}{x}+\sum_{n=0}^{\infty}(-1)^{n-1}\varsigma(n+1)x^n$ . We also remind

$$\chi_0(s)=\sum_{n=1}^{\infty}\psi(n)n^{-s}=\sum_{n=1}^{\infty}H_{n-1}n^{-s}-\gamma\varsigma(s)=\sum_{n=1}^{\infty}H_n n^{-s}-\gamma\varsigma(s)-\varsigma(s+1) \qquad (17)$$

whence $\chi_0(2n+1)=\sum_{k=1}^{\infty}\frac{H_k}{k^{2n+1}}-\gamma\varsigma(2n+1)-\varsigma(2n+2)$.

We have the following summation rule equating the coefficients standing in front of the $x^{2n-1}$, $n$=1, 2, 3…:

**Proposition 2.** For $n$=1, 2, 3…, we have

$$2n\sum_{k=1}^{\infty}\frac{H_k}{k^{2n+1}}=\sum_{l=1}^{2n-1}(-1)^{l-1}(2n-l)\varsigma(l+1)\varsigma(2n-l+1) \qquad (18).$$

Actually, already eq. (13) is sufficient to establish (18). This is *not* an Euler equality

$$\sum_{k=1}^{m-2}\varsigma(k+1)\varsigma(m-k)=(m+2)\varsigma(m+1)-2\sum_{n=1}^{\infty}\frac{H_n}{n^m} \qquad (19),$$

(also obtained by us in [1] working with the expansion $\psi(x)\psi(-x)+\frac{1}{x^2}-\gamma^2-\frac{\pi^2}{3}=2\sum_{n=1}^{\infty}\chi(2n+1)x^{2n}$ ) but can be shown to be equivalent to it when transformed applying known summation rules for Riemann zeta-functions. Moreover, the result of Proposition 2 can be considered as a proof of the following summation rule:

**Proposition 3.** For $n$=1, 2, 3…, we have

$$\frac{1}{n}\sum_{l=1}^{2n-1}(-1)^l(2n-l)\varsigma(l+1)\varsigma(2n-l+1)=\sum_{k=1}^{2n-1}\varsigma(k+1)\varsigma(2n+1-k)-(2n+3)\varsigma(2n+2) \qquad (20)$$

For example, the first two rules, following from (18) for $n$=1 and $n$=2 respectively, are $\sum_{n=1}^{\infty}\frac{H_n}{n^3}=\frac{\varsigma^2(2)}{2}=\frac{\pi^4}{72}$ and $2\sum_{n=1}^{\infty}\frac{H_n}{n^5}=2\varsigma(2)\varsigma(4)-\varsigma^2(3)=\frac{\pi^6}{270}-\varsigma^2(3)$, while from Euler summation rules (19) we have $\sum_{n=1}^{\infty}\frac{H_n}{n^3}=\frac{5}{2}\varsigma(4)-\frac{1}{2}\varsigma^2(2)=\frac{\pi^4}{72}$ and $2\sum_{n=1}^{\infty}\frac{H_n}{n^5}=7\varsigma(6)-2\varsigma(2)\varsigma(4)-\varsigma^2(3)=\frac{\pi^6}{270}-\varsigma^2(3)$, etc.

A bit more work is required for even powers of (15) or (16). We obtain first

$$\sum_{l=0}^{2n}(-1)^{l-1}(2n-l+1)\varsigma(l+1)\varsigma(2n-l+2)-(2n+3)\varsigma(2n+3)=-(2\chi_1(2n+1)-(2n+1)\chi_0(2n+2))$$
(21),

and then we have $\chi_0(2n+2)=\sum_{k=1}^{\infty}\frac{H_k}{k^{2n+2}}-\varsigma(1)\varsigma(2n+2)-\varsigma(2n+3)$

and $\chi_1(2n+1)=\varsigma(2)\varsigma(2n+1)-\sum_{k=1}^{\infty}\frac{H_k^{(2)}}{k^{2n+1}}+\varsigma(2n+3)$; see above. Thus we have proven the following Proposition.

**Proposition 4.** For $n$=1, 2, 3…, we have

$$\sum_{l=2}^{2n-1}(-1)^{l-1}(2n-l+1)\varsigma(l+1)\varsigma(2n-l+2)=(2n+1)\sum_{k=1}^{\infty}\frac{H_k}{k^{2n+2}}+2\sum_{k=1}^{\infty}\frac{H_k^{(2)}}{k^{2n+1}}-(2n+1)\varsigma(2)\varsigma(2n+1)$$
(22).

The first two summation rules following from (22) for $n$=1 and $n$=2 respectively are

$$3\sum_{k=1}^{\infty}\frac{H_k}{k^4}+2\sum_{k=1}^{\infty}\frac{H_k^{(2)}}{k^3}=3\varsigma(2)\varsigma(3) \qquad (23)$$

and

$$5\sum_{k=1}^{\infty}\frac{H_k}{k^6}+2\sum_{k=1}^{\infty}\frac{H_k^{(2)}}{k^5}=5\varsigma(2)\varsigma(5)-\varsigma(3)\varsigma(4) \qquad (24).$$

Equalities (23) and (24) have been tested numerically.

**Remark.** Terms $\sum_{k=1}^{\infty}\frac{H_k^{(2)}}{k^{2n+1}}$ can be isolated from (22) and then the values for the terms $\sum_{k=1}^{\infty}\frac{H_k}{k^{2n+2}}$ can be substituted using e.g. Euler summation rule to give an expression for the sum involving $H_n^{(2)}=\sum_{l=1}^{n}\frac{1}{l^2}$. For example, from (23) we have $2\sum_{k=1}^{\infty}\frac{H_k^{(2)}}{k^3}=3\varsigma(2)\varsigma(3)-3\sum_{k=1}^{\infty}\frac{H_k}{k^4}$, and thus using $\sum_{k=1}^{\infty}\frac{H_k}{k^4}=3\varsigma(5)-\varsigma(2)\varsigma(3)$:

$$2\sum_{k=1}^{\infty}\frac{H_k^{(2)}}{k^3}=6\varsigma(2)\varsigma(3)-9\varsigma(5) \qquad (25).$$

*1.3 Function $\psi^2(z)$*

Just for completeness, let us now briefly present the easiest example of the function $\psi^2(z)$. By squaring of eq. (4), we obtain $\psi^2(-(n+x)) = \frac{1}{x^2} + \frac{2(H_n - \gamma)}{x} + O(1)$. Thus for this function we have the poles of the second order, of the type *1/z*$^2$, and the poles of the first order with the residues $2(H_n - \gamma)$ at points *z*=0, -1, -2… The contour integral $\int_C z^{-s-1}\left(\psi^2(z) + \psi^2(-z) - \frac{2}{z^2} - 2\gamma^2 + 2\varsigma(2)\right)dz$ exists in the strip $0 < \mathrm{Re}\, s < 2$ and is equal to $2i\sin(\frac{s}{2}\pi)\int_0^\infty x^{-s-1}\left(\psi^2(ix) + \psi^2(-ix) + \frac{2}{x^2} - 2\gamma^2 + 2\varsigma(2)\right)dx$. From other side, by the Cauchy residues law, it is given by

$2\pi i\sum_{n=1}^\infty(\frac{s+1}{n^{s+2}} - 2\frac{H_n - \gamma}{n^{s+1}}) = 2\pi i((s+1)\varsigma(s+2) + 2\gamma\varsigma(s+1) - 2\sum_{n=1}^\infty \frac{H_n}{n^{s+1}})$. Repeating our standard consideration with inverse Mellin transform, we get

$$\psi^2(ix) + \psi^2(-ix) + \frac{2}{x^2} - 2\gamma^2 + 2\varsigma(2) = 2\sum_{n=1}^\infty(-1)^n((2n+1)\varsigma(2n+2) + 2\gamma\varsigma(2n+1) - 2\sum_{k=1}^\infty \frac{H_k}{k^{2n+1}})x^{2n}$$

and thus

$$\psi^2(x) + \psi^2(-x) - \frac{2}{x^2} - 2\gamma^2 + 2\varsigma(2) = 2\sum_{n=1}^\infty((2n+1)\varsigma(2n+2) + 2\gamma\varsigma(2n+1) - 2\sum_{k=1}^\infty \frac{H_k}{k^{2n+1}})x^{2n}.$$

Similarly, $\psi^2(x) - \psi^2(-x) + \frac{4\gamma}{x} = -2\sum_{n=1}^\infty(2n\varsigma(2n+1) + 2\gamma\varsigma(2n) - 2\sum_{k=1}^\infty \frac{H_k}{k^{2n}})x^{2n-1}$, and finally

$$\psi^2(x) - \frac{1}{x^2} - \gamma^2 + \varsigma(2) + \frac{2\gamma}{x} = \sum_{n=1}^\infty(-1)^n((n+1)\varsigma(n+2) + 2\gamma\varsigma(n+1) - 2\sum_{k=1}^\infty \frac{H_k}{k^{n+1}})x^n \qquad (26).$$

Comparing (26) with the square of (1) we get the following summation rule:

**Proposition 5.** For *n*= 2, 3…, we have

$$2\sum_{k=1}^\infty \frac{H_k}{k^n} = (n+2)\varsigma(n+1) - \sum_{k=1}^{n-2}\varsigma(n-k)\varsigma(k+1) \qquad (27).$$

Not surprisingly, here is nothing new: these are indeed Euler summation rules (19), see Proposition 2. The first of them (not appearing in our study above) is Euler equality $\sum_{k=1}^\infty \frac{H_k}{k^2} = 2\varsigma(3)$.

## 2. Conclusion

Following the program outlined in our publication [1], we established numerous summation rules involving the values of the Riemann zeta-function at integers and generalized harmonic numbers. Certainly, in all these rules the values of the Riemann

$$B_{2n} = \frac{(-1)^{n-1} 2(2n)!}{(2\pi)^{2n}} \varsigma(2n) \text{ [2].}$$

This line of research can be continued further by considering polygamma functions of larger order and their products, such as, say, $[\psi'(x)]^2 \pm [\psi'(-x)]^2$, or $\psi(x)\psi''(x) \pm \psi(-x)\psi''(-x)$ and so forth. This leads to the summation rules involving larger powers of harmonic numbers.

The present author indeed undertook some efforts in this direction, but obtained formulae are complicated and do not look very useful, so we do not present them here. Similarly, we continued also by considering e.g. the functions of the type $\psi(n+x) \pm \psi(n-x)$ or $\psi(n+x)\psi(n-x)$ - but again, obtained formulae are cumbersome and do not look truly interesting.